\documentclass[12pt]{amsart}  
\usepackage{amssymb,amsmath,amsthm,amscd}
\usepackage[all]{xy}

\usepackage{graphicx}
\usepackage{psfrag}

\numberwithin{equation}{section}

\newtheoremstyle{fancy1}{10pt}{10pt}{\itshape}{12pt}{\textsc\bgroup}{.\egroup}{8pt}{
}
\newtheoremstyle{fancy2}{10pt}{10pt}{}{12pt}{\itshape}{.}{8pt}{ }

\theoremstyle{fancy1}

\newtheorem{main}{Theorem}
\newtheorem*{main*}{Theorem}

\newtheorem*{cor*}{Corollary}

\newtheorem*{problem*}{Problem}

\theoremstyle{fancy2}

\newtheorem*{rem*}{Remark}

\newcommand{\cref}[1]{Corollary~\ref{#1}}

\newcommand{\CP}{\mathbb{C\mkern1mu P}}
\newcommand{\HP}{\mathbb{H\mkern1mu P}}

\newcommand{\Sph}{\mathbb{S}}

\newcommand{\R}{{\mathbb{R}}}

\newcommand{\QH}{{\mathbb{H}}}

\newcommand{\G}{\ensuremath{\operatorname{G}}}

\newcommand{\SO}{\ensuremath{\operatorname{SO}}}

\newcommand{\Sp}{\ensuremath{\operatorname{Sp}}}
\newcommand{\U}{\ensuremath{\operatorname{U}}}
\newcommand{\SU}{\ensuremath{\operatorname{SU}}}
\newcommand{\Spin}{\ensuremath{\operatorname{Spin}}}

\newcommand{\T}{\ensuremath{\operatorname{T}}}
\renewcommand{\S}{\ensuremath{\operatorname{S}}}

\newcommand{\No}{\ensuremath{\operatorname{N}}}

\newcommand{\fg}{{\mathfrak{g}}}
\newcommand{\fk}{{\mathfrak{k}}}

\newcommand{\fm}{{\mathfrak{m}}}

\newcommand{\fsu}{{\mathfrak{su}}}

\newcommand{\fsp}{{\mathfrak{sp}}}

\newcommand{\fspin}{{\mathfrak{spin}}}
\newcommand{\fspi}{{\mathfrak{spin}}}

\def\con#1=#2(#3){#1 \equiv #2 \bmod{#3}}

\newcommand{\tr}{\ensuremath{\operatorname{tr}}}

\renewcommand{\Im}{\ensuremath{\operatorname{Im}}}

\renewcommand{\sec}{\ensuremath{\operatorname{sec}}}

\DeclareMathOperator{\rpart}{Re}

\newcommand{\no}{\noindent}

\begin{document}

\title{Positively curved homogeneous metrics on spheres}

\author{Luigi Verdiani}
\address{Universit\'a di Firenze\\ Via S. Marta 3\\50139 Firenze, Italy}
\email{luigi.verdiani@unifi.it}

\author{Wolfgang Ziller}
\address{University of Pennsylvania\\
        Philadelphia, PA 19104}
\email{wziller@math.upenn.edu}

\thanks{The first named author was supported by GNSAGA. The last named author was supported in part by the
Francis J. Carey Term Chair, the Clay Institute and a grant from the
National Science Foundation.}

\maketitle

The study of manifolds with positive sectional curvature has a long
history and can be considered as the beginning of  global Riemannian
geometry. Nevertheless there are few general theorems concerning
this class of manifolds which do not require further geometric
assumptions. Apart from obstructions which already hold for
manifolds with non-negative sectional curvature, there are only
restrictions on the fundamental group due to the classical theorems
of Bonnet-Meyers and Synge.

There are also very few known examples. They all arise as quotients
of a compact Lie group, endowed with a left invariant metric, by a
subgroup of isometries acting freely. They consist, apart from the
rank one symmetric spaces, of certain homogeneous spaces in
dimensions $6,7,12,13$ and $24$ due to Berger \cite{Be}, Wallach
\cite{Wa}, and Aloff-Wallach \cite{AW}, and of biquotients in
dimensions $6,7$ and $13$ due to Eschenburg \cite{E1},\cite{E2} and
Bazaikin \cite{Ba}.

\smallskip

The homogeneous spaces which admit homogeneous metrics with positive
sectional curvature have been classified (\cite{Wa,BB}). A natural
further question is the behavior of the pinching constants, i.e. the
quotient of the minimum by the maximum of the sectional curvature.
Of particular interest is whether one can determine in each case the
homogeneous metric with the best pinching constant. This necessarily
involves a classification of all positively curved homogeneous
metrics on each space, and the determination of all 2-planes with
minimal and maximal sectional curvature. This is a notoriously
difficult algebraic problem and has been solved (numerically in some
cases) in a remarkable paper by P\"uttmann \cite{Pu} for all
homogeneous spaces not diffeomorphic to a sphere (see also
\cite{El,Gr,He,Hua,Va} for previous particular cases).

\smallskip

In this paper we will consider the only remaining case, that of
homogeneous metrics on spheres, and will  classify the set of
positively  curved ones. Apart from the natural challenge that this
question poses, we will see that the set of curvature tensors
involved is surprisingly complicated, and the methods developed for
solving the problem may be of interest in other situations as well.

Homogeneous metrics on spheres were classified in \cite{Zi} and can
be described geometrically in terms of the Hopf fibrations:

$$\Sph^1\to \Sph^{2n+1}\to \CP^n\; , \;
\Sph^3\to \Sph^{4n+3}\to \HP^n\; ,\; \Sph^7\to \Sph^{15}\to
\Sph^8.$$ By scaling the round sphere metric on the total space by a
constant $t$ in the direction of the fibers, we obtain a one
parameter family of metrics $g_t$. For this family of metrics we
will show:

\begin{main}
The homogeneous metrics $g_t$ have positive sectional curvature if
and only if $\; 0<t<4/3$ and their pinching constants are given in
\eqref{max}.
\end{main}

This result is well known in the case where the fiber is one
dimensional since in that case the curvature operator is diagonal
and the maximum and minimum of the sectional curvature is hence
obtained at what we call {\it natural} two-planes. They consist of
vertical, horizontal, and vertizontal two planes, depending on
wether  the two plane is spanned by two vertical vectors, two
horizontal ones, or one vertical and one horizontal. If the fiber
dimension is 3 or 7,  the curvature operator is not diagonal anymore
and if $1<t<4/3$ the maximum, and if $4/5<t<1$ the minimum of the
sectional curvatures is  attained at two planes whose projection
onto the horizontal and the vertical space are both two dimensional.

\smallskip

The only remaining homogeneous metrics are given (up to scaling) by
the family $g_{(t_1,t_2,t_3)}$ on $\Sph^{4n+3}$ where we modify the
round sphere metric with an arbitrary left invariant metric on the
3-sphere fiber. Up to isometry, one can assume the metric is
diagonal with respect to a fixed basis and $t_i$ are the lengths
squared of the basis vectors.

To describe the positive curvature conditions, let
$$V_i= (t_j^2+t_k^2-3t_i^2+2
t_it_j+2 t_it_k-2 t_jt_k)/{t_i} \; \text{ and } \; H_i=4-3t_i,$$
with $(i,j,k)$ a cyclic permutation of $(1,2,3)$.  We will show:

\begin{main}
The homogeneous metrics $g_{(t_1,t_2,t_3)}$ on $\Sph^{4n+3}$ have
positive sectional curvature if and only if
$$ V_i>0 \quad , \quad H_i>0 \quad \text{ and }
\quad 3|t_jt_k-t_j-t_k+t_i| < t_jt_k +\sqrt{H_i V_i} $$ with
$(i,j,k)$ a cyclic permutation of $(1,2,3)$.
\end{main}

The cone condition $V_i>0$  ensures that the totally geodesic fibers
have positive sectional curvature and $H_i>0$ guarantees that all
horizontal curvatures are positive. The third condition is more
subtle and again  due to two-planes with a two dimensional
projection onto the vertical and onto the horizontal space. This
extra condition cuts out a very thin slice from the star shaped
region $V_i>0\ , \ H_i>0$, see Figures 3-6.

\smallskip

In the proof of Theorem A we determine all critical points of the
sectional curvature function. For the 3-dimensional family in
Theorem B, this direct approach turns out to be intractable. The
method used in \cite{Pu} (which is due to Thorpe \cite{Th})  works
by modifying the curvature operator with invariant 4-forms. Although
this approach was very successful for the homogeneous metrics
studied  in \cite{Pu}, it turns out that it gives little information
for the metrics $g_{(t_1,t_2,t_3)}$. Instead, we will determine
those metrics where the minimum of the  sectional curvature is 0, by
explicitly solving for the two planes that achieve this minimum.

\smallskip

This work was completed while the second author was visiting IMPA
and he would like to thank the Institute for its hospitality.

\section{Preliminaries}

If $G$ is a compact Lie group acting transitively on
$\Sph^{n-1}\subset \R^n$ then  the action of $G$ is equivalent to a
linear action of a subgroup of $\SO(n)$ (see \cite{MS}, \cite{Bo}).
The set of $\G$ invariant metrics was determined in \cite{Zi} and
can be described, up to scaling, as follows:

\begin{itemize}
\item $\G=\SO(n) , \Spin(7)$ or $\G_2$. The isotropy representation
is irreducible and hence all metrics have constant curvature.
 \item
$\G=\U(n)$ or $\SU(n)$. The isotropy representation consists of two
irreducible summands of dimension $1$ and $2n-2$ and the $\G$
invariant metrics are the round sphere metric scaled in the
direction of the Hopf fibration with one dimensional fibers.
\item $\G=\Sp(n)\Sp(1)$ or $\G=\Spin(9)$. The isotropy representation consists of two
irreducible summands of dimension $3$ and $4n-4$ respectively $7$
and $8$, and the $\G$ invariant metrics are the round sphere metric
scaled in the direction of the Hopf fibration with three
respectively seven dimensional fibers.
\item $\G=\Sp(n)$. The isotropy representation consists of a trivial
$3$-dimensional and a $4n-4$ dimensional irreducible summand and the
$\G$ invariant metrics are the round sphere modified on the fibers
of the Hopf fibration with three dimensional fibers by an arbitrary
left invariant metric. Up to isometries the metric on the fiber can
be assumed to be diagonal.
\item $\G=\Sp(n)\U(1)$. This is the subclass of metrics in the previous case
where the metric on the fiber has two diagonal entries of the same
length.
\end{itemize}

\bigskip

Since the metrics in the case of $\G=\U(n)$ have a diagonal
curvature operator, it remains to consider two cases:

\begin{itemize}
\item[(a)] A one parameter family of metrics where the Hopf
fibration with $3$ or $7$-dimensional fibers is scaled by $t$ in the
direction of the fibers of a Hopf fibration.
\item[(b)] A three parameter family of metrics where the Hopf
fibration with $3$-dimensional fibers is modified in the direction
of the fibers by a left invariant metric.
\end{itemize}

\no We will discuss case (a) in Section 2, and case (b) in Section
3.

\bigskip

We first describe the metrics in the case of
$\Sph^{4n-1}=\Sp(n)/\Sp(n-1)$. On the Lie algebra of $\Sp(n)$ we
have the biinvariant inner product $g_0( A,B) = -\frac 1 2
\rpart(\tr AB)$. We then consider the orthogonal splitting
$\fsp(n)=\fsp(n-1)\oplus\fm$ and identify $\fm$ with the tangent
space of $\Sp(n)/\Sp(n-1)$ at the identity coset. A  basis of $\fm$,
orthonormal in $g_0$, is given by $X_{r}\, , r=1\cdots 3$ and $
U_{rs}\ , r=1\cdots 4\, ,s=1\cdots n-1$. Here $X_{r}$ are the
matrices which have as its only non-zero entries,
$\sqrt{2}i,\sqrt{2}j$ respectively $\sqrt{2}k$ in  row one and
column one. $U_{1s}$ has as its only non-zero entries a $1$ in row 1
and column $s+1$ and a $-1$ in row $s+1$ and column $1$. $U_{2s}$
has as its only non-zero entries an $i$ in row 1 and column $s+1$
and in row $s+1$ and column $1$ as well. Similarly $U_{3s}$ and
$U_{4s}$ with $i$ replaced by $j$ or $k$ respectively. In this
basis, the metric $g_0$ induced on $\Sph^{4n-1}$ has constant
curvature 1. The vectors $X_r$ are tangent to the fibers of the Hopf
fibration and $U_{rs}$ is a basis of the horizontal space. The
isotropy group $\Sp(n)$ acts trivially on the vertical space and
irreducibly on the horizontal space. We can therefore scale the
metric such that it is equal to $g_0$ on the horizontal space. On
the vertical space the inner product is arbitrary, but as was
observed in \cite{Zi}, we can use the gauge group
$\No(\Sp(n))/\Sp(n)=\Sp(1)$ to make the metric on the vertical space
diagonal. In other words, the vectors $X_i$ are orthogonal and have
length squared $t_i\, , i=1\cdots 3$. We will denote this metric by
$g_{(t_1,t_2,t_3)}$. If $t_1=t_2$, the metrics are invariant under
the bigger group $\Sp(n)\U(1)$ and if $t_1=t_2=1$ they are invariant
under $\U(2n)$. Finally, if $t_1=t_2=t_3:=t$, the metric is
invariant under $\Sp(n)\Sp(1)$ and will be denoted by $g_t$.

\bigskip

In the case of $\Sph^{15}=\Spin(9)/\Spin(7)$, the basis is more
complicated and can be chosen as follows. Let $E_{ij}$ be the
standard basis of $\fspin(9)$, and let
\begin{align*}
\fk_1& = \text{ span } \{ E_{24}+E_{68}, E_{28}+E_{46}, E_{26}-E_{48}\} \\
\fk_2& =\text{ span }\{ E_{23}+E_{67}, E_{27}+E_{36}, E_{34}+E_{78},
E_{38}+E_{47}, E_{37}-E_{48}\} \\
 \fk_3& =\text{ span }\{ E_{27}-E_{45},
E_{23}+E_{58}, E_{24}-E_{57}, E_{28}+E_{35}, E_{56}-E_{78},
2\,E_{25}-E_{38}+E_{47}\} \\
\fk_4& =\text{ span }\{ E_{12}+E_{56}, E_{16}+E_{25}, E_{13}+E_{57},
E_{17}+E_{35}, E_{14}+E_{58}, E_{18}+E_{45}, E_{15}-E_{48}\}\\
\fm_1& =\text{ span }\{ X_1=E_{15}+E_{26}+E_{37}+E_{48},
X_2=E_{17}+E_{28}-E_{35}-E_{46},\\
& \hspace{49pt} X_3 =E_{13}-E_{24}-E_{57}+E_{68},
X_4=E_{16}-E_{25}-E_{38}+E_{47},\\
& \hspace{49pt} X_5=E_{18}-E_{27}+E_{36}-E_{45},
 X_6=E_{12}+E_{34}-E_{56}-E_{78},\\
& \hspace{49pt} X_7=E_{14}+E_{23}-E_{58}-E_{67}\}\\
\fm_2& =\text{ span }\{ U_1=E_{19}, U_2=E_{29}, U_3=E_{39},
U_4=E_{49},
U_5=E_{59}, U_6=E_{69},\\
& \hspace{49pt} U_7 =E_{79}, U_8=E_{89}\}.
\end{align*}

The basis is chosen such that $\fk_1\simeq \fsu(2)$,
$\fsu(2)+\fk_2\simeq \fsu(3)$, $\fsu(3)+\fk_3\simeq \fg_2$, and
$\fg_2+\fk_4\simeq \fspi(7)$. The orthogonal complement to
$\fspin(7)$ is $\fm_1\oplus \fm_2$ and is identified with the
tangent space. The isotropy group leaves $\fm_i$ invariant and acts
irreducibly on them. Thus we can scale any $\Spin(7)$ invariant
inner product such that $U_i$ is an orthonormal basis of $\fm_2$,
$\fm_1$ and $\fm_2$ are orthogonal, and the vectors $V_i\in\fm_1$
are orthogonal with length squared equal to $t$. This metric will be
denoted by $g_t$.

\smallskip

Throughout the paper, $i,j,k$ in any formula will always denote  a
cyclic permutation of $(1,2,3)$.

\section{Pinching constants for $g_t$}

We now study the metrics $g_t$ and compute their pinching. We start
with the case where the metric is invariant under $G=\Sp(n)\Sp(1)$
and use the basis $X_r ,U_{r,s}$ of the tangent space described in
Section 1. They are all orthogonal to each other, $\vert X_i\vert^2
=t$ and the remaining vectors have unit length. The isotropy group
$\Sp(n-1)\Sp(1)$ acts on the tangent space as follows:
$(A,q)(y,v)=(qyq^{-1},Avq^{-1})$ where $y\in \Im \QH$ and $v\in
\QH^n$. Modulo this isotropy action we can thus assume that the
2-plane is spanned by $X,Y$ with
\begin{align}
X&=a_1\, X_1+a_2\, X_2+a_3 X_3+a_4 U_{11}\\
\notag Y&=b_1\, X_1+b_2\, X_2+\, b_3\, X_3+b_4\, U_{11}+b_5\,
U_{21}+b_6\,U_{12}.
\end{align}

A computation shows that (see also Section 3):
\begin{align*}
\langle R(X,Y)X,Y\rangle&=t \{(a_2b_3-a_3b_2)^2+  (a_3b_1-a_1b_3)^2+
(a_1b_2-a_2b_1)^2\}\notag \\
&\quad +(4-3t)\, a_4^2b_5^2+6t(t-1)\, (a_2b_3-a_3b_2)\, a_4b_5\\
&\quad + t^2\{ (b_4a_1-b_1a_4)^2+  (b_4a_2-b_2a_4)^2 +
(b_4a_3-b_3a_4)^2\}\notag \\
&\quad + t^2b_5^2(a_1^2+  a_2^2 + a_3^2)+b_6^2
(t^2(a_1^2+a_2^2+a_3^2)+a_4^2) .\notag
\end{align*}

The sectional curvatures of the natural 2-planes are $1/t$ for
vertical 2-planes, $t$ for vertizontal 2-planes and they lie in
$[4-3t,1]$ for horizontal 2-planes. Here we use the convention that
in an interval $[a,b]$ we allow $a\le b$ or $b\le a$. If  we want
the sectional curvatures to be positive, we need to assume that
 $0<t<4/3$, and we will see that this is indeed also sufficient.
 In addition we want to compute the pinching and hence we minimize and maximize
 all sectional curvatures. Our strategy
will be to see if all  curvatures lie in the interval given by the
values at natural 2-planes. We will see that this is not always the
case and will determine the 2-planes that actually achieve the
maximum and minimum.

\smallskip

  If we set
$A=(a_1,a_2,a_3)$ and $B=(b_1,b_2,b_3)$,
 the sectional curvature is  given by
\begin{align*}
\sec(X,Y)&\{(t\vert A\vert^2+a_4^2)(t\vert
B\vert^2+b_4^2+b_5^2)-(t(A\cdot B)+a_4b_4)^2+b_6^2( t\vert
A\vert^2+a_4^2)\} \\
&=t \vert A\times B\vert^2 +(4-3t)\, a_4^2b_5^2
-6t(t-1)\, ((A\times B)\cdot e_1)\, a_4b_5+\\
&\quad + t^2\{ b_4^2\vert A\vert^2+a_4^2\vert B\vert^2
-2a_4b_4(A\cdot B) + b_5^2\vert A\vert^2\} +b_6^2 (t^2\vert
A\vert^2+a_4^2) .\notag
\end{align*}

In the following we will use the observation that a function of the
type $\frac{ax^2+b}{cx^2+d}$, where $a,b,c,d$ are constants, is
either decreasing or increasing on $[0, \infty]$. We can apply this
for the variables $b_4$ and $b_6$ if we subtract a multiple of $X$
from $Y$ in order to make $A\cdot B=0$, which we will assume from
now on. In that case, $\sec\to t$ when $b_4\to \infty$ and $\sec\to
\frac{t^2|A|^2+a_4^2}{t|A|^2+a_4^2}\subset [1,t]$ when $b_6\to
\infty$, and we can hence assume that $b_4=b_6=0$. Finally, we can
make  $(A\times B)\cdot e_1$ be equal to any value in the interval
$[-\vert A\vert\vert B\vert,\vert A\vert \vert B\vert \ ]$ without
changing $\vert A\vert$ or $ \vert B\vert$.
 Since $a_4$ can take on negative values, we
can thus assume that $(A\times B)\cdot e_1=\vert A\vert\vert B\vert
$ in order to determine the maximum and minimum of the sectional
curvatures.

If we set $x=\vert A\vert,\, y=\vert B\vert,\, r=a_4 ,\, s= b_5$ we
hence have to determine the critical points of
$$
F=\frac{t x^2y^2 +(4-3t)\, r^2s^2 +6t(t-1)\, xyrs+
t^2(x^2s^2+y^2r^2)}{(tx^2+r^2)(ty^2+s^2)}.
$$

Using the above observation, one sees again that if any of the 4
variable vanish, the values of $F$ lie between the curvatures of
natural 2-planes. We can thus normalize the vectors such that
$x=y=1$ and one easily sees that the only non-zero critical points
of the remaining function of $r$ and $s$ are, besides $(r,s)=(0,0)$,
given by
$$
r=-s=\pm \sqrt{1-2t} \qquad \text{ with }\qquad  F=4(1-t)
$$
and
$$
r=s=\pm \sqrt{(4t+1)/t}\quad  \text{ with } \qquad
F=\alpha:=\frac{16t^2-8t+4}{11t+1}.
$$
The value of $F$ for the first critical point lies in between the
curvatures of natural 2-planes for $0<t<1/2$. But for the second
critical point, $\alpha>t$ for $1<t<4/3$  and  $\alpha<t$ for
$4/5<t<1$. We thus obtain for the maximum of the sectional curvature
\begin{equation*}
 \max \sec =
\begin{cases} \ 1/t & \text{ if \ \ $0<t\le 1/3$,}
\\
\ 4-3t &\text{ if \ \ $1/3\le t\le 1$,}
\\ \ \alpha & \text{ if \ \ $1\le t\le 4/3$,}
\end{cases}
\end{equation*}
and for the minimum
\begin{equation*}
 \min \sec =
\begin{cases} \ t & \text{ if \ \  $0<t\le 4/5$,}
\\
\ \alpha &\text{ if \ \ $4/5\le t\le 1$,}
\\ \  4-3t & \text{ if \ \  $1\le t\le 4/3$.}
\end{cases}
\end{equation*}
Hence  the pinching constants, i.e.,  $\delta=\min \sec/\max \sec$,
are given by
\begin{equation}
\label{max} \hspace{20pt} \text{ pinching } \delta =
\begin{cases} \  t^2 & \text{ if \ \  $0<t\le 1/3$,}
\\
\ t/(4-3t) &\text{ if \ \  $1/3\le t\le 4/5$,}
\\
\ \alpha/(4-3t) &\text{ if \ \  $4/5\le t\le 1$,}
\\  \ (4-3t)/ \alpha & \text{ if \ \  $1\le t\le 4/3$.}
\end{cases}
\end{equation}

Figure 1 shows a graph of the pinching number on the left. For $4/5<
t<4/3$ the pinching is slightly less than what one would obtain by
using only natural 2-planes. The difference is shown in the picture
on the right and is less than $0.008$.

\psfrag{1}{$\scriptstyle 1 $\ }

\psfrag{1.01}{$\scriptstyle 1 $\ }

\psfrag{t}{$\scriptstyle t $\ }

\psfrag{0.33}{$\scriptstyle \frac 1 3$}

\psfrag{1.33}{$\scriptstyle \frac 4 3$}

\psfrag{0.501}{$\scriptstyle \frac 1 2$}

\psfrag{0.111}{$\scriptstyle \frac 1 9$}

\psfrag{0.00350001}{$\scriptstyle 0.003$}

\psfrag{0.00800001}{$\scriptstyle 0.008$}
\psfrag{delta_t}{$\scriptstyle \delta_t$}
\begin{figure}[h]
\begin{center}
\psfrag{0.8}{$\scriptstyle \frac 4 5$}
\includegraphics[width=2.8in,height=2.8in]{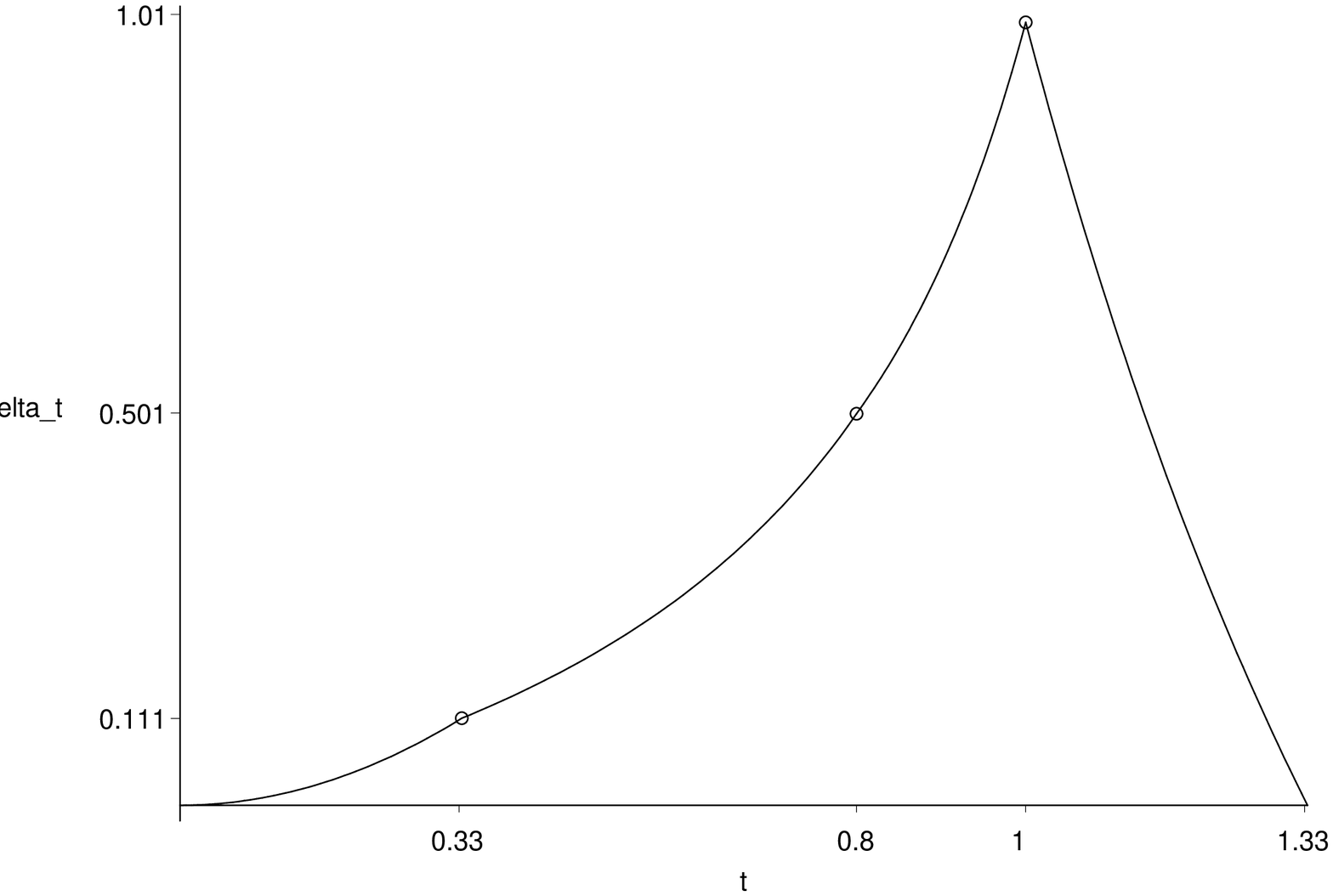}\qquad \quad
\includegraphics[width=2.8in,height=2.8in]{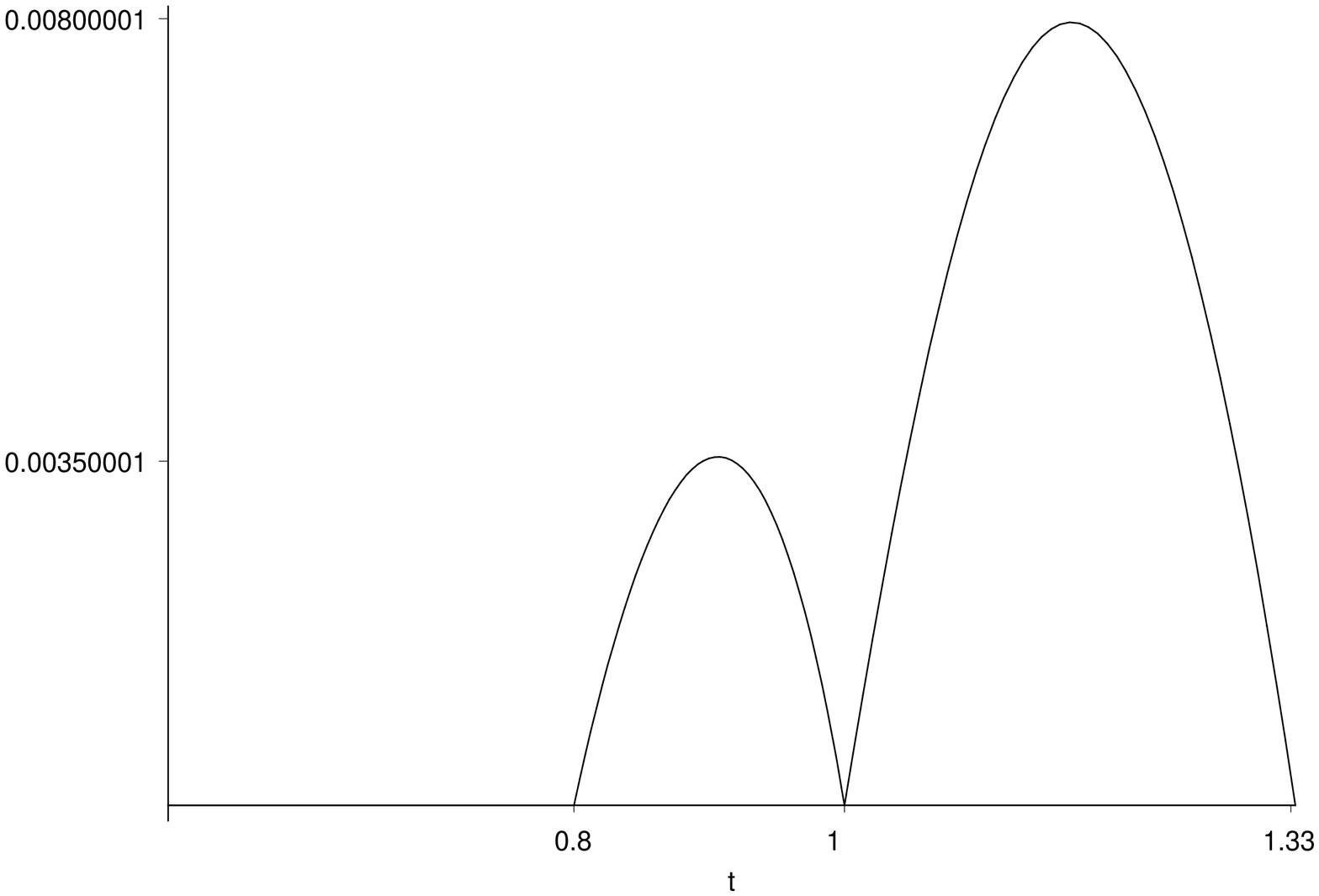}
\end{center}
\caption{Pinching for $g_t$ and difference with natural 2-planes.}
\end{figure}

\bigskip

We now claim that the same conclusions hold for
$\Sph^{15}=\Spin(9)/\Spin(7)$. We will use the basis for the tangent
space described in Section 1. Let $X=A_1+B_1$ ,$Y=A_2+B_2$ be a
basis of a $2$-plane in $\fm$ with $A_i\in \fm_1$, $B_i\in \fm_2$.
The isotropy group $\Spin(7)$ acts via the standard representation
on $\fm_1$ and via the spin representation on $\fm_2$. Thus, since
$\Spin(7)/\G_2=\Sph^7$ and $\G_2/\SU(3)=\Sph^6$ , we may assume
$B_1\in \text{span}(U_1)$, $B_2\in \text{span}(U_1,U_5)$. Since the
isotropy group $\SU(3)$ fixes $X_1$,   we may hence assume $A_1\in
\text{span}(X_1,X_2)$, although for symmetry reasons we will choose
$A_1\in \text{span}(X_1,X_2,X_3)$. Since $\SU(2)$ fixes
$\text{span}(X_1,X_2,X_3)$, we may choose $A_2\in
\text{span}(X_1,X_2,X_3,X_4)$.

Let $X=a_1\,X_1+a_2\,X_2+a_3\,X_3+a_4\,U_1$ and
$Y=b_1\,X_1+b_2\,X_2+b_3\,X_3+x\,X_4+b_4\,U_1+b_5\,U_5$. If we set
$A=(a_1,a_2,a_3)$ and $B=(b_1,b_2,b_3)$ as before,  the sectional
curvature is  given by
\begin{align*}
\sec(X,Y)&\{(t\vert A\vert^2+a_4^2)(t\vert
B\vert^2+b_4^2+b_5^2)-(t(A\cdot B)+a_4b_4)^2+x^2( t\vert
A\vert^2+a_4^2)\} \\
&=t \vert A\times B\vert^2 +(4-3t)\, a_4^2b_5^2
+6t(t-1)\, ((A\times B)\cdot e_1)\, a_4b_5+\\
&\quad + t^2\{ b_4^2\vert A\vert^2+a_4^2\vert B\vert^2
-2a_4b_4(A\cdot B) + b_5^2\vert A\vert^2\} +x^2(t\vert
A\vert^2+t^2a_4^2),\notag
\end{align*}
and the argument proceeds as before.  In particular, for both cases
of $g_t$, maximum and minimum are assumed at 2-planes tangent to a
totally geodesic 7-sphere.

\section{Positive curvature on $S^{4n-1}$}

In this section we will study the curvature tensor of the metrics
$g_{(t_1,t_2,t_3)}$ on $S^{4n-1}$ described in Section 1 and use the
basis $X_r,U_{rs}$ for the tangent space. Modulo the action of the
isotropy subgroup $\Sp(n-1)$, it is enough to consider 2-planes
spanned by vectors of the form
\begin{align}
X&=a_1\, X_1+a_2\, X_2+a_3\, X_3+a_4\, U_{11}\\
\notag
 Y&=b_1\, X_1+b_2\,
 X_2\, +b_3\, X_3\, +b_4\, U_{11}+b_5\, U_{21}+b_6\, U_{31}+b_7\,
U_{41}+b_8\,U_{12}.
\end{align}
Besides the quantities $V_i$ and $ H_i$ defined in the introduction,
we also set
$$
L_i=6(t_jt_k-t_j-t_k+t_i).
$$

 Using any one of the curvature formulas for a homogeneous space
 (see, e.g., \cite{Pu}),
one shows that
\begin{align*} &R(X_i,X_j,X_i,X_j)=V_k\; ,\;
R(U_{1p},U_{\gamma p},U_{1p},U_{\gamma p})
= R(U_{\alpha p},U_{\beta p},U_{\alpha p},U_{\beta p})=H_\gamma \\
 &R(U_{rp},U_{rq},U_{rp},U_{rq})=R(U_{rp},U_{sq},U_{rp},U_{sq})=1\ ,
 \
 R(U_{1p},U_{\alpha p},U_{\beta p},U_{\gamma p})
=\sigma (t_\beta+t_\gamma-2t_\alpha )\\
&R(U_{1p},U_{\alpha p},U_{1q},U_{\alpha q})= R(U_{\beta p},U_{\gamma
p},U_{\beta q},U_{\gamma q})=2(1-t_\alpha) \ , \  R(U_{1p},U_{\alpha
p},U_{\beta
q},U_{\gamma q})=\sigma (2-t_\alpha)\\
&R(U_{1p},U_{1q},U_{\alpha p},U_{\alpha q})=R(U_{1p},U_{\alpha
q},U_{1q},U_{\alpha p}) =R(U_{\beta p},U_{\beta q},U_{\gamma
p},U_{\gamma q})=1-t_\alpha\\
&R(U_{\alpha p},U_{\beta q},U_{\alpha q},U_{\beta p})= (1-t_\gamma)
\ , \ R(U_{1p},U_{\alpha q},U_{\beta p},U_{\gamma q}
=-\sigma (1-t_\beta)\\
&R(U_{1p},X_\alpha,U_{1p},X_\alpha)= R(U_{\beta p},X_\alpha,U_{\beta
p},X_\alpha)
=t_\alpha^2\\
&R(U_{\alpha p},U_{\beta p},X_\alpha,X_\beta) =2R(U_{\alpha
p},X_\alpha,U_{\beta p},X_\beta)
=2R(U_{\alpha p},X_\beta,X_\alpha,U_{\beta p})=- {L_\gamma}/3\\
&R(U_{1p},U_{\alpha p},X_\beta,X_\gamma)=
-2R(U_{1p},X_\gamma,U_{\alpha p},X_\beta) = -\sigma  {L_\alpha}/3,
\end{align*}
where $p\neq q\in \{1,(n-1)\}$, $r\neq s\in \{1,2,3,4\}$ and
$(\alpha,\beta,\gamma)$  a permutation, with sign $\sigma$, of
$(1,2,3)$. This easily implies that

\begin{align}\label{curv}
\langle R(X,Y)X,Y\rangle&=V_1 (a_2b_3-a_3b_2)^2+ V_2
(a_3b_1-a_1b_3)^2+ V_3
(a_1b_2-a_2b_1)^2\notag \\
&\quad +L_1\, (a_2b_3-a_3b_2)\, a_4b_5+ L_2\, (a_3b_1-a_1b_3)\,
a_4b_6+
L_3\, (a_1b_2-a_2b_1)\, a_4b_7\notag \\
 &\quad +H_1\, a_4^2b_5^2+ H_2\,
a_4^2b_6^2+
H_3\, a_4^2b_7^2 \\
&\quad + t_1^2 (b_4a_1-b_1a_4)^2+ t_2^2 (b_4a_2-b_2a_4)^2 +t_3^2
(b_4a_3-b_3a_4)^2\notag \\
&\quad + (b_5^2+b_6^2+b_7^2)(t_1^2 a_1^2+ t_2^2 a_2^2 +t_3^3
a_3^2)+  b_8^2(t_1^2 a_1^2+t_2^2 a_2^2+t_3^2 a_3^2+a_4^2)
 .\notag
\end{align}

 We now consider the set of homogeneous metrics
$$\mathcal P=\{g_{(t_1,t_2,t_3)}\mid t_i>0\, , V_i>0\, , H_i>0 \}
$$
 where the vertical and the horizontal curvatures are already
positive. Notice that \eqref{curv} implies that then all vertizontal
curvature are automatically positive as well.
  We want to check if these conditions
are also sufficient by looking for metrics in $\mathcal P$ where the
minimum of the sectional curvature is $0$. We thus assume from now
that $g=g_{(t_1,t_2,t_3)}\in \mathcal P$ is such that for any choice
of the parameters $a_i , b_i$ we have $R(X,Y,X,Y)\geq 0$, with
equality for some value of $a_i$ and $b_i$. In the following, we
will assume that $X,Y$ span such a zero curvature 2-plane and  our
strategy will be to try to solve for $a_i , b_i$ in terms of $t_i$.

We can clearly assume that $b_8=0$ since otherwise the two plane
spanned by $\bar{X} =X , \bar{Y} = Y-b_8 \, U_{21}$ would have
negative sectional curvature contradicting our assumption. Thus we
can restrict ourselves to metrics on $\Sph^7$.

 The expression of the curvature in \eqref{curv} suggests that one should consider
 vectors $X,Y$ orthonormal with respect to a modified inner product. We will assume in the
  following that $X$ and $Y$ satisfy the  conditions
\begin{equation}
\label{restr} \left\{
\begin{array}{l}
t_1^2 a_1^2+t_2^2 a_2^2+t_3^2 a_3^2+a_4^2=1\medskip\\
t_1^2 b_1^2+t_2^2 b_2^2+t_3^2 b_3^2+b_4^2+b_5^2+b_6^2+b_7^2=1\medskip\\
t_1^2 a_1b_1+t_2^2 a_2b_2+t_3^2 a_3b_3+a_4b_4=0.
\end{array}
\right.
\end{equation}
For later purposes we  observe that, under these conditions, we can
assume that $0<a_4<1$. Indeed, if $a_4=0$, either one of $b_i\ ,
i=4, \dots , 7$ is non-zero or the 2-plane is vertical. In either
case, \eqref{curv} implies that all curvatures are positive since
$g\in \mathcal P$. Similarly, if $a_4=1$, the curvatures are
positive since $a_i=0 , \ i=1,2, 3$ and thus  one of $b_i, \ i=1 , 2
, 3$ or $ i=5, 6 , 7$ is non-zero.

Using \eqref{restr}, we may  rewrite the last two lines in
(\ref{curv}) as follows:
\begin{align*}
t_1^2 &(b_4a_1-b_1a_4)^2+ t_2^2 (b_4a_2-b_2a_4)^2 +t_3^2
(b_4a_3-b_3a_4)^2 \\
& \hspace{15pt} + (b_5^2+b_6^2+b_7^2)(t_1^2 a_1^2+ t_2^2 a_2^2
+t_3^3 a_3^2)\\
& =  b_4^2 (t_1^2 a_1^2+t_2^2 a_2^2+t_3^2 a_3^2)+
a_4^2(t_1^2b_1^2+t_2^2b_2^2+t_3^2b_3^2)\notag \\
  & \hspace{15pt} -2a_4b_4(t_1^2a_1b_1+t_2^2a_2b_2+t_3^2a_3b_3)
+(b_5^2+b_6^2+b_7^2) (t_1^2 a_1^2+t_2^2 a_2^2+t_3^2 a_3^2)\\
& = (b_4^2+b_5^2+b_6^2+b_7^2) (1-a_4^2)+
a_4^2(1-b_4^2-b_5^2-b_6^2-b_7^2)-
2a_4b_4(-a_4b_4)\notag\\
& = b_4^2+a_4^2+(b_5^2+b_6^2+b_7^2) (1-2 a_4^2).\\
\end{align*}
We now claim that under the assumptions \eqref{restr}, the value of
$b_4$ must be $0$. Indeed, if we let $\overline X=X$ and $\overline
Y=Y-b_4V_1$ we obtain
\begin{align*}
\langle R(X,Y)X,Y\rangle-\langle R(\overline X,\overline Y)\overline
X,\overline Y\rangle=b_4^2(1-a_4^2)+2a_4^2b_4^2=b_4^2(1+a_4^2)
\end{align*}
and hence $\langle R(\overline X,\overline Y)\overline X,\overline
Y\rangle<0$, contradicting our assumption.

\smallskip

Thus we have to minimize the following function, where we have used
the further abbreviation $A_{ij}=a_ib_j-a_jb_i$
\begin{align}\label{function}
F(a_i,b_i)=&V_1A_{23}^2+V_2A_{31}^2+V_3A_{12}^2+  L_1 A_{23}\,
a_4b_5+L_2A_{31}\, a_4b_6+L_3A_{12}\, a_4b_7 \notag \\
&+H_1\, a_4^2b_5^2+H_2\, a_4^2b_6^2+H_3\, a_4^2b_7^2+
a_4^2+(b_5^2+b_6^2+b_7^2) (1-2 a_4^2),
\end{align}
subject to the restrictions \eqref{restr}. We use  Lagrange
multipliers to determine the minimizing $0$ curvature 2-planes. The
Lagrange multiplier equations with respect to $a_1,a_2,a_3$ are
given by
\begin{align*}
&2 V_3 b_2A_{12}-2 V_2 b_3 A_{31}+L_3 a_4b_2b_7-L_2a_4b_3b_6+2 t_1^2
\lambda_1 a_1+t_1^2 \lambda_3 b_1=0 \notag \\
 &2 V_1 b_3A_{23}-2 V_3 b_1
A_{12}+L_1 a_4b_3b_5-L_3a_4b_1b_7+2 t_2^2 \lambda_1 a_2+t_2^2
\lambda_3 b_2=0  \\
 &2 V_2 b_1A_{31}-2 V_1 b_2 A_{23}+L_2
a_4b_1b_6-L_1a_4b_2b_5+2 t_3^2 \lambda_1 a_3+t_3^2 \lambda_3 b_3=0.
\notag
\end{align*}
By setting $S_i=2\,V_i\,A_{jk}+L_i\,a_4\,b_{4+i}$, we may rewrite
the  above equations as
\begin{align*}
&S_3\,b_2-S_2\,b_3+2 t_1^2 \lambda_1 a_1+t_1^2 \lambda_3 b_1=0\\
&S_1\,b_3-S_3\,b_1+2 t_2^2 \lambda_1 a_2+t_2^2 \lambda_3 b_2=0\\
&S_2\,b_1-S_1\,b_2+2 t_3^2 \lambda_1 a_3+t_3^2 \lambda_3 b_3=0,
\end{align*}
that is,
$$-(S_1,S_2,S_3)\times (b_1,b_2,b_3)+2\,
\lambda_1\,(t_1^2\,a_1,t_2^2\,a_2,t_3^2\,a_3)+\lambda_3\,
(t_1^2\,b_1,t_2^2\,b_2,t_3^2\,b_3)=(0,0,0).$$ Since $(b_1,b_2,b_3)$
is orthogonal to $(t_1^2\,a_1,t_2^2\,a_2,t_3^2\,a_3)$ this  implies
that $(b_1,b_2,b_3)$ is orthogonal also to
$\lambda_3\,(t_1^2\,b_1,t_2^2\,b_2,t_3^2\,b_3)$. If $\lambda_3\neq
0$ it follows that $b_1=b_2=b_3=0$. But in this case \eqref{curv}
implies that the sectional curvature is positive, contradicting our
assumption. Hence $\lambda_3=0$ and
$(t_1^2\,a_1,t_2^2\,a_2,t_3^2\,a_3)$ is orthogonal to
$(S_1,S_2,S_3)$.

In a similar way, using the Lagrange multipliers for $b_1$, $b_2$,
$b_3$  we also have the orthogonality between
 $(t_1^2\,b_1,t_2^2\,b_2,t_3^2\,b_3)$ and $(S_1,S_2,S_3)$. Hence there exists a $\lambda$ (we use
 $-2\lambda$ for convenience)
  with
\begin{align*}\label{lambda}
(S_1,S_2,S_3)&= -2\lambda (t_1^2\,a_1,t_2^2\,a_2,t_3^2\,a_3)\times
(t_1^2\,b_1,t_2^2\,b_2,t_3^2\,b_3)   \\
&= -2\lambda(A_{23}t_2^2 t_3^2 ,\ A_{31}t_3^2t_1^2 ,\
A_{12}t_1^2t_2^2). \notag
\end{align*}
We can thus rewrite the Lagrange multiplier equations for
$a_1,a_2,a_3,b_1,b_2,b_3$ as
\begin{equation}
\label{Eq1} \left\{
\begin{array}{l}
L_1 a_4b_5=-2(\lambda t_2^2t_3^2+V_1)A_{23}\medskip\\
L_2 a_4b_6=-2(\lambda t_3^2t_1^2+V_2)A_{31}\medskip\\
L_3 a_4b_7=-2(\lambda t_1^2t_2^2+V_3)A_{12}.\medskip
\end{array}
\right.
\end{equation}

We can now  determine the multipliers. We claim they are given by
\begin{equation}
\label{lambda2} \lambda_1=\lambda (1-b_5^2-b_6^2-b_7^2) \quad ,\quad
\lambda_2=\lambda (1-a_4^2)\quad , \quad \lambda_3=0
\end{equation}
To see this, say for $\lambda_1$,  set
$$v=(t_1^2 a_1,t_2^2 a_2,t_3^2 a_3)\ ,
w=(t_1^2 b_1,t_2^2 b_2,t_3^2 b_3)\, , \ z=(b_1,b_2,b_3).$$ We thus
have
$$2\lambda_1 v=(S_1,S_2,S_3)\times z=-2\lambda (v\times w)\times z.$$
Under the restrictions \eqref{restr}, it follows that $v\cdot z=0$
and $ w\cdot z = 1- b_5^2-b_6^2-b_7^2$. Using the identity $(v\times
w)\times z = (w\cdot z)v-(v\cdot z)w$ the claim follows.

\bigskip

Consider now the multiplier equations with respect to $b_5,b_6,b_7$:
\begin{align*}
2b_5(1-2a_4^2)+2b_5 H_1 a_4^2+L_1 A_{23} a_4+2\lambda_2 b_5=0\\
2b_6(1-2a_4^2)+2b_6 H_2 a_4^2+L_2 A_{31} a_4+2\lambda_2 b_6=0\\
2b_7(1-2a_4^2)+2b_7 H_3 a_4^2+L_3 A_{12} a_4+2\lambda_2 b_7=0
\end{align*}
Using the abbreviation
$$h=\frac {1-2\,a_4^2+\lambda_2}{a_4}$$
we can rewrite the equations as
\begin{equation}
\label{Eq2} \left\{
\begin{array}{l}
H_1 a_4 b_5=-b_5 h-  \frac{L_1}{2}A_{23}\medskip\\
H_2 a_4 b_6=-b_6 h-  \frac{L_2}{2}A_{31}\medskip\\
H_3 a_4 b_7=-b_7 h-  \frac{L_3}{2}A_{12}\medskip.

\end{array}
\right.
\end{equation}

The Lagrange multiplier equation for the variable $a_4$ is given by
$$
\begin{aligned}
& 2 H_1 a_4 b_5^2+2 H_2 a_4 b_6^2+2 H_3 a_4 b_7^2 +L_1 A_{23}
b_5+L_2A_{31}b_6+L_3A_{12}b_7 \\ & + 2 a_4-4
a_4(b_5^2+b_6^2+b_7^2)+2 \lambda_1 a_4=0.
\end{aligned}
$$
Using \eqref{Eq2} it follows that
$$
 -2(b_5^2+b_6^2+b_7^2)\,h +2 a_4-4 a_4(b_5^2+b_6^2+b_7^2)+2
\lambda_1 a_4=0.
$$
If in turn we replace the value of $\lambda_1 , h$ and $\lambda_2$
in this expression, and multiplying by $a_4$,  the multiplier
equation for $a_4$ becomes
\begin{equation}
\label{Eq3} (a_4^2-b_5^2-b_6^2-b_7^2)(\lambda+1)=0.
\end{equation}

We have as a last condition that the  sectional curvature of the
2-plane is $0$. Using \eqref{function}- \eqref{Eq2} we get
\begin{align*}
0 & =V_1 A_{23}^2+V_2A_{31}^2+V_3A_{12}^2+ \frac {L_1}{2} a_4b_5
A_{23}+\frac {L_2}{2} a_4b_6
A_{31}+\frac {L_3}{2} a_4b_7 A_{12}\\
& \hspace{15pt} -h a_4(b_5^2+b_6^2+b_7^2)
+a_4^2+(b_5^2+b_6^2+b_7^2)(1-2a_4^2) \\ & =V_1
A_{23}^2+V_2A_{31}^2+V_3A_{12}^2 -(\lambda t_2^2t_3^2+V_1)A_{23}^2
-(\lambda t_3^2t_1^2+V_2)A_{31}^2 -
(\lambda t_1^2t_2^2+V_3)A_{12}^2\\
& \hspace{15pt} +a_4^2+(b_5^2+b_6^2+b_7^2)(1-2a_4^2-h a_4)\\
& = - \lambda (A_{23}^2t_2^2 t_3^2 + A_{31}^2t_3^2t_1^2 +
A_{12}^2t_1^2t_2^2)+a_4^2-(b_5^2+b_6^2+b_7^2)\lambda_2\\
& = -\lambda (1-b_5^2-b_6^2-b_7^2)(1-a_4^2) + a_4^2-\lambda
(b_5^2+b_6^2+b_7^2) (1-a_4^2)\\
& = a_4^2-\lambda (1-a_4^2).
\end{align*}
In the above computation we have used the consequence of
\eqref{restr} that the vector with components $t_ia_i , i=1,2 , 3$
and the one with components $t_ib_i , i=1,2 , 3$ are orthogonal to
each other and hence their cross product has length squared
$(1-b_5^2-b_6^2-b_7^2)(1-a_4^2)$. It thus follows that
$$\lambda=\frac{a_4^2}{1-a_4^2}$$
which in turn implies  $$h=\frac {1-a_4^2}{a_4} .$$ Furthermore,
comparing with (\ref{Eq3}), we obtain
$$a_4^2=b_5^2+b_6^2+b_7^2.$$

\smallskip

Altogether, we see that we can find a minimizing $2$-plane with zero
sectional curvature if and only if we can solve the following
system:

\begin{equation}
\label{system} \left\{
\begin{array}{l}
L_i a_4b_{4+i}=-2\left(t_j^2t_k^2 \frac{a_4^2}{1-a_4^2}+V_i\right)A_{jk}\medskip\\
H_i a_4b_{4+i}=-\frac{1-a_4^2}{a_4^2}a_4b_{4+i}-\frac {L_i}{2} A_{jk}\medskip\\
a_4^2=b_5^2+b_6^2+b_7^2\medskip
\end{array}
\right. \qquad i=1,2,3
\end{equation}
under the restrictions \eqref{restr}. Since we  assumed that $a_4>
0$, there must be one $b_i , i>4$ which does not vanish. If
$b_{4+i}\neq 0$, we can solve the first equation for
$\frac{A_{jk}}{a_4b_{4+i}}$, substitute into the second, and obtain
a quadratic equation in $Z=\frac{a_4^2}{1-a_4^2}$:

\begin{equation}
\label{quad} t_j^2t_k^2 H_i Z^2+\left(t_j^2t_k^2 +H_i V_i-\frac
{L_i^2}{4}\right)Z+V_i=0.
\end{equation}

In this equation, the leading coefficient and the last term are
positive by assumption. Since we also assumed that $a_4<1$, we need
a positive solution,  which implies in particular that
$E_i:=t_j^2t_k^2 +H_i V_i-\frac {L_i^2}{4} <  0$. We now claim that
there can be at most one $i$ with $E_i < 0$. For this purpose, we
will show that  $E_1> 0$ if $0<t_1\leq t_2\le 4/3$. It then follows,
by symmetry, that the same is true for $t_1\leq t_3$ and hence $E_1$
can be negative only in the region where $t_1>\max\{t_2,t_3\}$.
Similarly for $E_2$ and $E_3$  in the regions $t_2>\max\{t_1,t_3\}$
and $t_3>\max\{t_1,t_2\}$ respectively, which will finish our claim.

The fact that $E_1>0$ for $0<t_1\leq t_2\le 4/3$ can be proved by
somewhat tedious estimates. We indicate the reason via a picture. We
factor $E_1$ as follows:

\begin{align*}
t_1E_1=At_3^2+Bt_3+C&=
(-12t_1+4-8t_1t_2^2+18t_1t_2)t_3^2\\
& +(-12t_1t_2-8t_2+12t_1^2+18t_1t_2^2
+8t_1-18t_1^2t_2)t_3\\
&-12t_1^2+4t_2^2+12t_1^2t_2-12t_1t_2^2+8t_1t_2.
\end{align*}

One can now draw a graph of the coefficients $A$ and $C$ to check
that  in the region $0<t_1\leq t_2\le 4/3$ the function $A$ is
positive and $C$ is non-negative. Hence $t_1E_1$ is an upward
pointing parabola in the $t_3$ variable which thus has as its
minimum $C-\frac{B^2}{4A}$. Figure 2 shows where this minimum is
positive in the white area in the right hand side picture. The left
hand side picture shows in the black area where the coefficient $B$
is negative. Altogether, this implies that $E_1$ is positive.

\begin{figure}[h]
\begin{center}
\includegraphics[width=2.5in,height=2.5in]{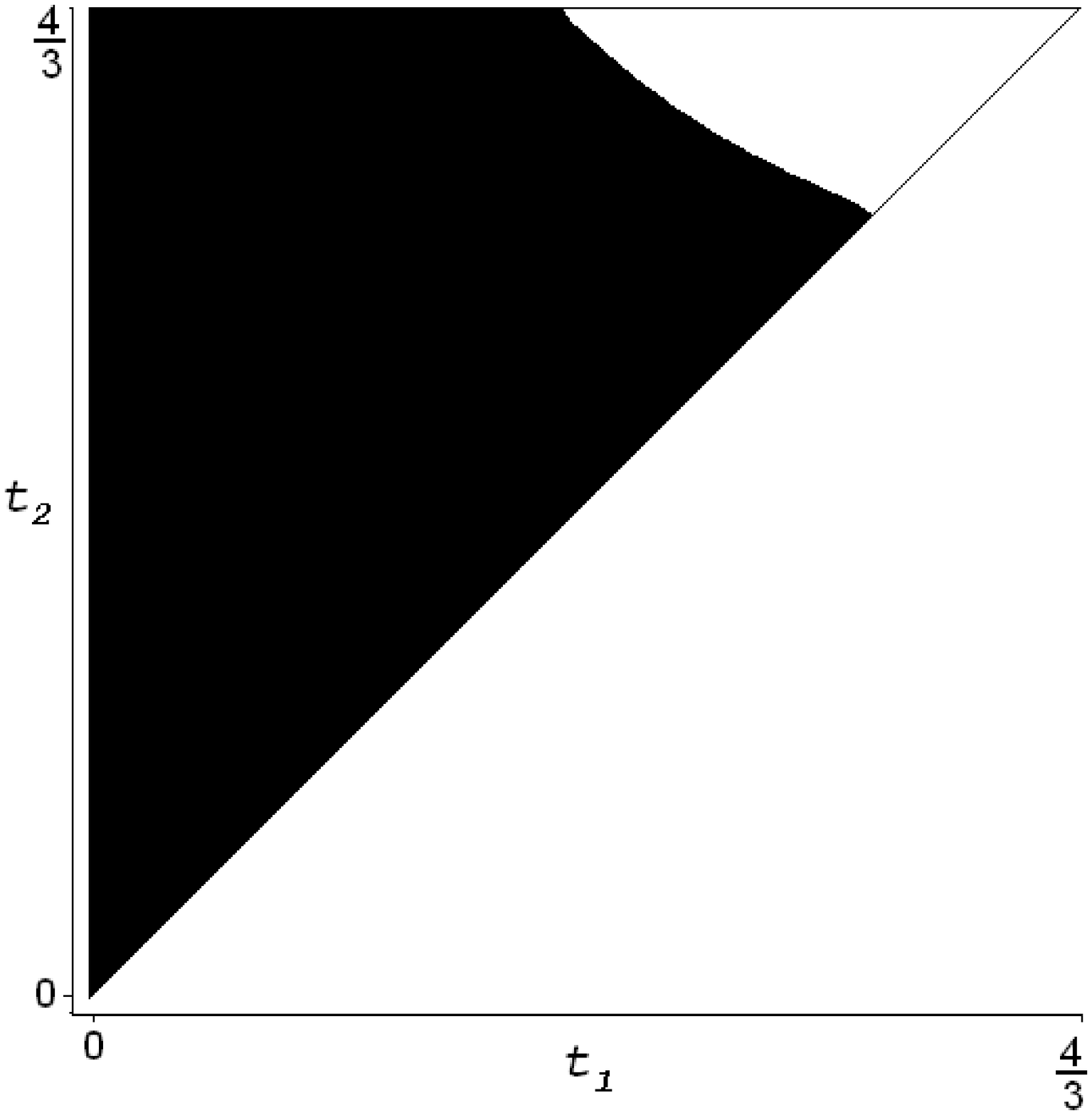}\qquad \quad
\includegraphics[width=2.5in,height=2.5in]{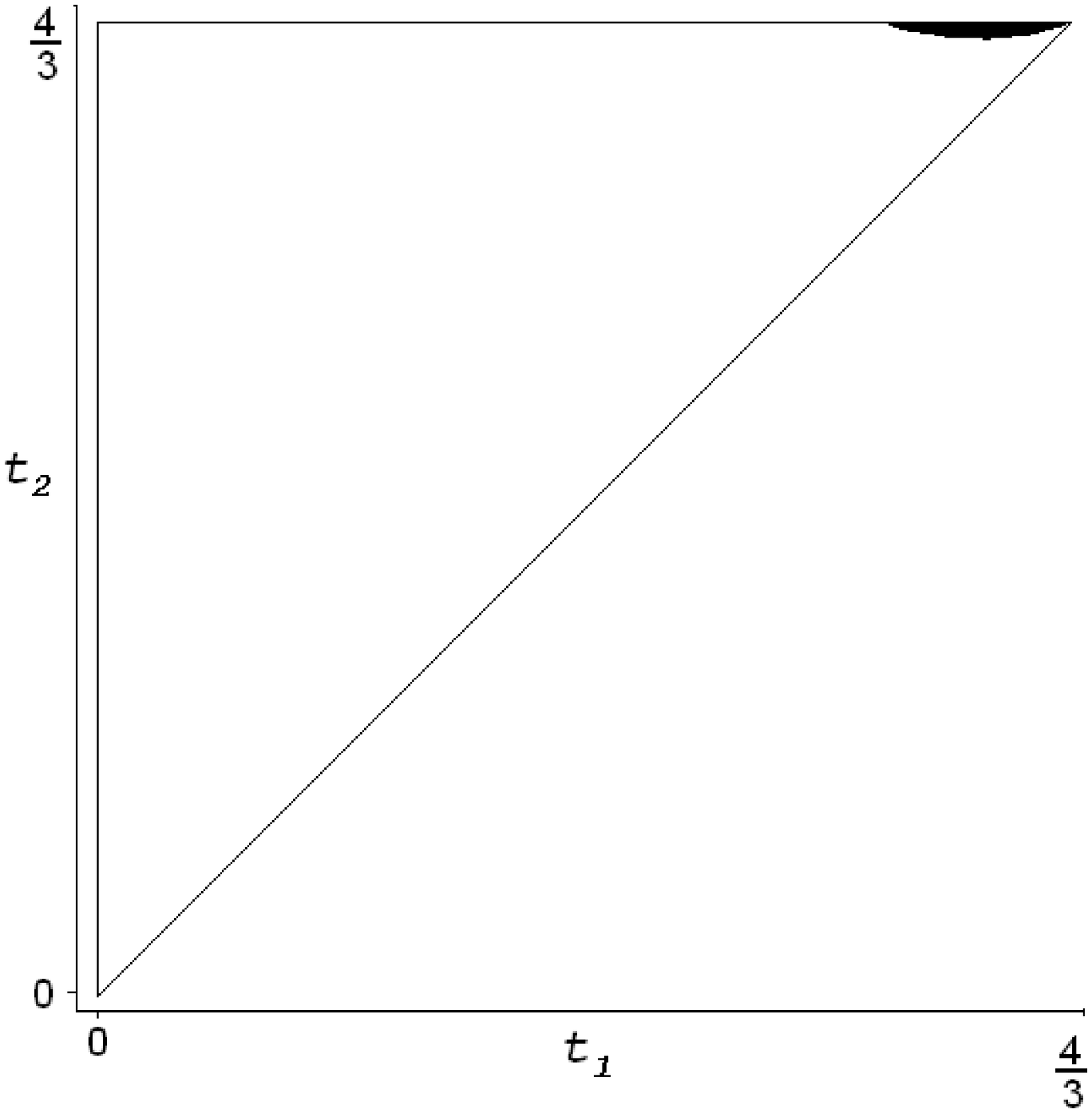}
\end{center}
\caption{}
\end{figure}

Thus we have shown that only one of the equations \eqref{quad} can
have a positive solution. Assume that this is the equation for
$i=1$, in other words $b_5\ne 0$. We then necessarily have
$b_6=b_7=0$ and hence $b_5^2=a_4^2$. Furthermore, \eqref{system}
with $i=2,3$ can only hold if $A_{13}=a_1b_3-a_3b_1=0$ and
$A_{12}=a_1b_2-a_2b_1=0$ and hence $a_1=b_1=0$. Thus the
normalization \eqref{restr} takes on the form
$$ t_2^2a_2^2+t_3^2a_3^2=t_2^2b_2^2+t_3^2b_3^2=1-a_4^2,
\qquad t_2^2a_2b_2+t_3^2a_3b_3=0,$$ which implies that
$(t_2b_2,t_3b_3)=\pm (t_3a_3,-t_2a_2)$ and hence
$A_{23}=a_2b_3-a_3b_2=\mp \frac{1-a_4^2}{t_2t_3}$.

 Instead of solving the system \eqref{system} directly, we can
  now substitute the values of
$a_i$ and $b_i$ into the sectional curvature formula $(3.4)$ and
obtain:
\begin{equation}\label{finsec}
\frac{t_2^2t_3^2}{(1-a_4^2)^2}F=t_2^2t_3^2H_1Z^2
 + 2t_2t_3 (\pm \frac {L_1}2+t_2t_3) Z+V_1.
\end{equation}
Since we require that $F\ge 0$,  a positive solution $Z$ to $F=0$
corresponds to a double root which exists  if and only if
$$|L_1|=2(t_2t_3+\sqrt{H_1V_1}).$$

Altogether we have shown that a metric in $\mathcal P$  with
non-negative sectional curvature admits a 0 curvature plane if and
only if $|L_i| = 2(t_jt_k +\sqrt{H_i V_i})$ for one $i$. Finally,
observe that \eqref{finsec} also implies that a metric with $|L_i|
>2(t_jt_k +\sqrt{H_i V_i})$ has 2-planes with negative sectional
curvature. Hence a metric in $\mathcal P$ with $\sec>0$ is
characterized by $|L_i| < 2(t_jt_k +\sqrt{H_i V_i})$, which proves
Theorem B.

\bigskip

\begin{rem*}
If we have a metric with non-negative, but not positive curvature,
the proof also classifies the set of $2$-planes with $0$ curvature.
One easily sees that the $0$-curvature $2$-planes are all vertical
if $V_i=0$ and for $H_i=0$ they are all horizontal. For the
remaining metrics in the boundary of $\mathcal P$ the $0$-curvature
$2$-planes are described as follows. As we saw, only one of the
quadratic equations  \eqref{quad} has a positive solution. Assuming
that this is the case for  $i=1$, the 0 curvature  planes  consist
of two disjoint circles where each $2$-plane is spanned by:
$$X=\frac{\cos(\theta)}{t_2}X_2 + \frac{\sin(\theta)}{t_3}X_3 +
\sqrt{\frac{Z}{1+Z}}U_1\, ,\quad Y=\frac{-\sin(\theta)}{t_2}X_2 +
\frac{\cos(\theta)}{t_3}X_3 \pm \sqrt{\frac{Z}{1+Z}}U_2
$$
for some $\theta$, and  where $Z$ is a positive solution of
\eqref{quad}.
\end{rem*}

We illustrate  the set of positively curved metrics in Theorem B
with some pictures. Figure 3 shows the cone of metrics where
horizontal and vertizontal curvatures are positive and Figure 4 a
cross section in the plane $t_1+t_2+t_3=1$. The cone touches the
coordinate planes in the diagonals $t_i=t_j\;, t_k=0$. Thus in this
case $t_k$ can go to $0$, corresponding to a Cheeger deformation in
the direction of the extra circle of isometries, i.e. contracting in
the orbits of the circle action. Such Cheeger deformations preserve
positive curvature.

Figure 5 shows the surfaces $V_1=0$ and in dark the extra small
slice cut out by the additional inequality in Theorem B. Figure 6
shows the difference in the $t_3$ coordinates of the two surfaces as
a function of $t_1,t_2$. The difference is at most $0.0085$. The curve
that separates the two surfaces is given by $4t_1t_2 = 4t_1+4t_2-3$.
A typical metric with some negative sectional curvatures in between
these two surfaces is given by $(t_1,t_2,t_3)=(0.25,0.25,0.33)$. The
2-plane spanned by $(a_1,a_2,a_3,a_4)=(4,\; 0,\; 0,\; \frac 7 {10})$
and $(b_1,b_2,b_3,b_4,b_5,b_6,b_7)=$ \newline $(0,\; 4,\; 0,\; 0,\;
0,\; 0,\; \frac 7 {10})$ then has  $\sec(X,Y)=-0.69$.

Figure 7 shows a cross section of the set of positively curved
metrics by the plane $t_1=t_2$ which corresponds to $\Sp(n)\U(1)$
invariant metrics. The dotted curve goes from $(0,0)$ to $(\frac 1 2
,\frac 2 3 )$ and is given by the equation $t_3=\dfrac{
t_1(4t_1^3-12t_1^2-4+9t_1)}{3(2t_1-2t_1^2-1)}$. It is slightly less
than $t_3=\frac 4 3 t_1$ which corresponds to $V_3=0$. Figure 8
again shows the difference between these two curves. A dot in the
pictures represents the biinvariant metric.

\psfrag{t1}{$\scriptstyle t_1$}\psfrag{t2}{$\scriptstyle
t_2$}\psfrag{t3}{$\scriptstyle t_3$} \psfrag{1.33333}{$\scriptstyle
 \frac 4 3$}
\begin{figure}[h]
\begin{center}
\includegraphics[width=4in,height=4in]{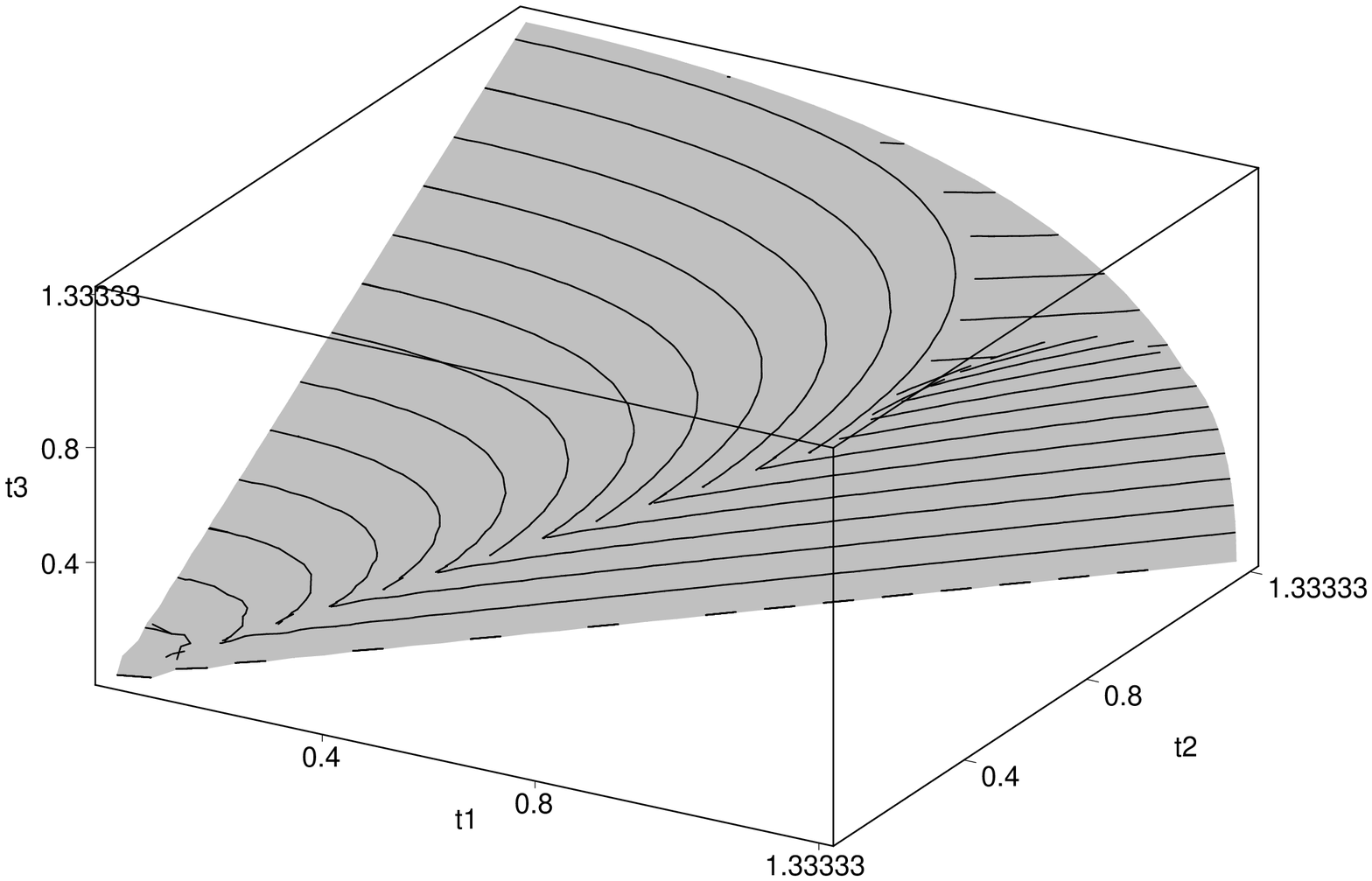}
\end{center}
\caption{Positive cone $\mathcal P=\{g_{t_1,t_2,t_3} \mid V_i>0 , H_i > 0\} $.}
\end{figure}

\bigskip

\begin{figure}[h]
\begin{center}
\includegraphics[width=3in,height=3.3in,angle=-90]{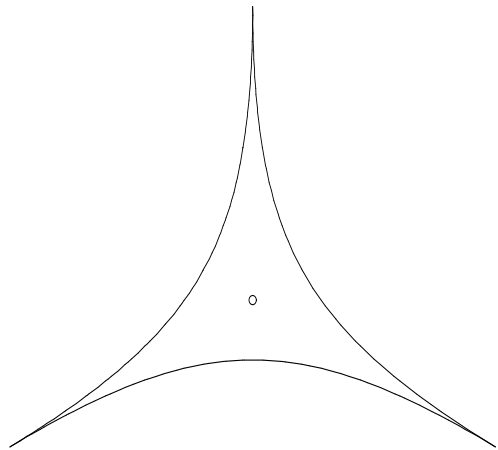}
\end{center}
\caption{Cross section of $\mathcal P$ in 2-plane $t_1+t_2+t_3=1$.}
\end{figure}

\bigskip

\begin{figure}[h]
\begin{center}
\includegraphics[width=4in,height=4in,angle=-90]{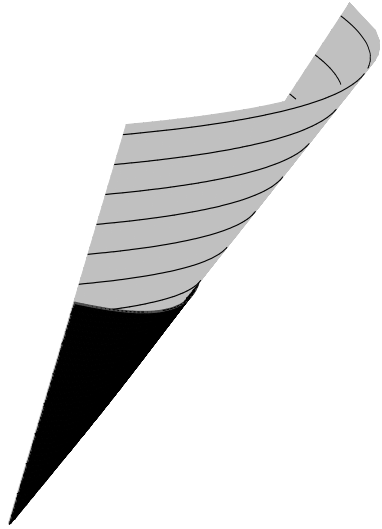}
\end{center}
\caption{One face of modified $\mathcal P$.}
\end{figure}

\bigskip

\begin{figure}[h]
\begin{center}
\includegraphics[width=3.3in,height=3.3in,angle=-90]{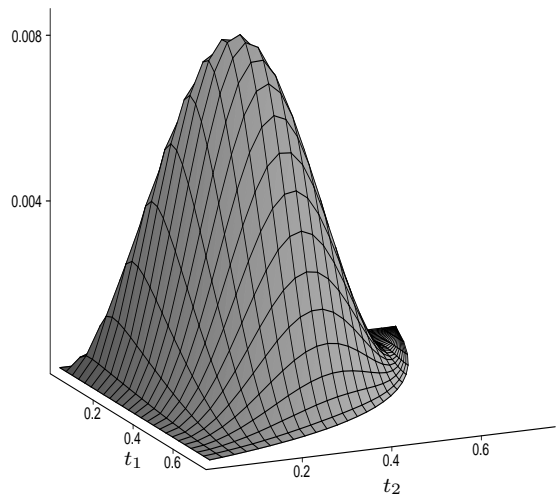}
\end{center}
\caption{Difference of the two surfaces.}\label{disjoint}
\end{figure}

\bigskip

\psfrag{1.33}{$\scriptstyle \frac 4 3$} \psfrag{0.66}{$\scriptstyle
\frac 2 3$}\psfrag{0.5}{$\scriptstyle \frac 1 2$}

\begin{figure}[h]
\begin{center}
\includegraphics[width=4in,height=4in,angle=-90]{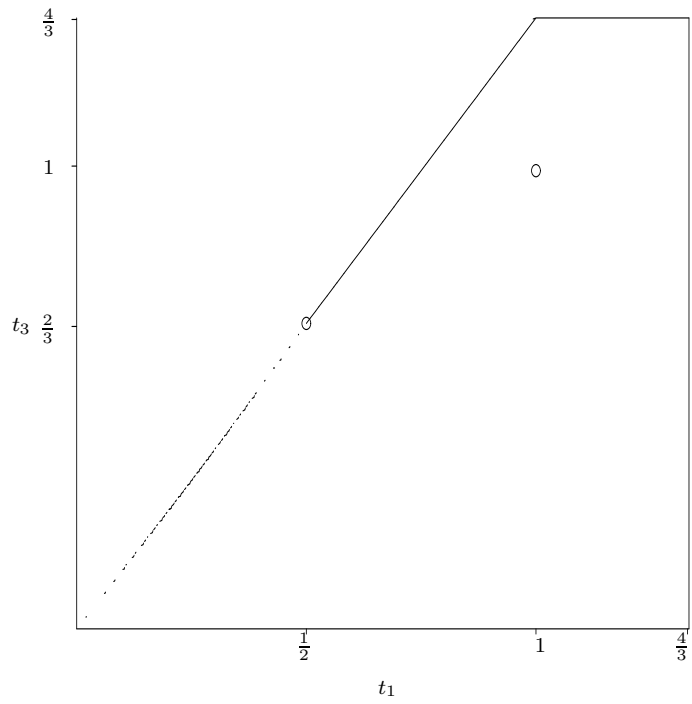}
\end{center}
\caption{Slice with $t_1=t_2$.}
\end{figure}

\bigskip

\begin{figure}[h]
\begin{center}
\includegraphics[width=3.1in,height=3.1in,angle=-90]{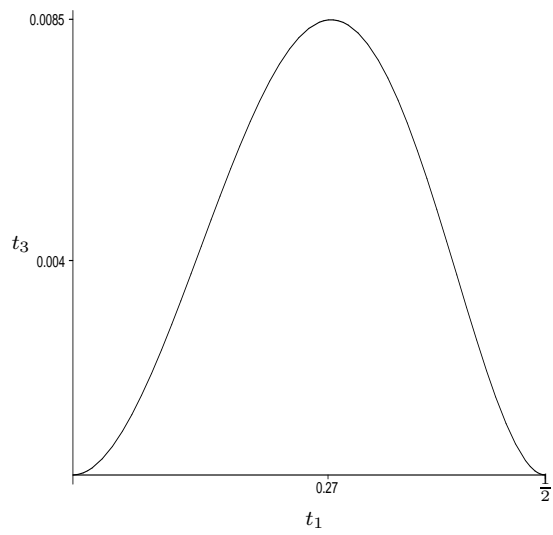}
\end{center}
\caption{Difference in height.}
\end{figure}

\end{document}